\documentclass[reqno,11pt]{amsproc}
\usepackage{amssymb}
\usepackage{graphicx}
\usepackage{amscd}
\usepackage{amsmath}
 \theoremstyle{plain}

\newtheorem{definition}{Definition}

\newtheorem{theorem}{Theorem}

\numberwithin{equation}{section}

\begin{document}

\title[Dirac system on a finite tree.]{Inverse dynamic
and spectral problems for the one-dimensional Dirac system on a
finite tree.}

\author{Alexander Mikhaylov} % required
\address{St. Petersburg   Department   of   V.A. Steklov    Institute   of   Mathematics
of   the   Russian   Academy   of   Sciences, 7, Fontanka, 191023
St. Petersburg, Russia and Saint Petersburg State University,
St.Petersburg State University, 7/9 Universitetskaya nab., St.
Petersburg, 199034 Russia.} \email{mikhaylov@pdmi.ras.ru}

\author{Victor Mikhaylov} % required
\address{St.Petersburg   Department   of   V.A.Steklov    Institute   of   Mathematics
of   the   Russian   Academy   of   Sciences, 7, Fontanka, 191023
St. Petersburg, Russia and Saint Petersburg State University,
St.Petersburg State University, 7/9 Universitetskaya nab., St.
Petersburg, 199034 Russia.} \email{ftvsm78@gmail.com}

\author{Gulden Murzabekova}%
\address{S. Seifullin Kazakh Agrotechnical University, 62 Pobeda Prospect, Astana, 010011, Kazakhstan}%
\email{guldenmur07@gmail.com}%

\keywords{Inverse problem, planar tree, Dirac system,
Titchmarsh-Weyl matrix, leaf peeling method}
%\date{September, 2017}

\maketitle

\begin{abstract}
We consider inverse dynamic and spectral problems for the one
dimensional Dirac system on a finite tree. Our aim will be to
recover the topology of a tree (lengths and connectivity of edges)
as well as matrix potentials on each edge. As inverse data we use
the Weyl-Titchmarsh matrix function or the dynamic response
operator.

\end{abstract}

\section{Introduction}

Let $\Omega$ be a finite connected compact graph without cycles (a
tree). The graph consists of edges $E=\{e_1,\ldots,e_{N}\}$
connected at the vertices $V=\{v_1\ldots,v_{N+1}\}$. Every edge
$e_j\in E$ is identified with an interval $(0,l_j)$ of the real
line. The edges are connected at the vertices $v_j$ which can be
considered as equivalence classes of the edge end points, we write
$e\sim v$ if the vertex $v$ is a boundary of the edge $e$. The
boundary $\widetilde\Gamma=\{v_1,\ldots,v_m\}$ of $\Omega$ is a
set of vertices having multiplicity one (the exterior nodes). In
what follows we assume that one boundary node (say $v_m$) is
clamped, i.e. zero Dirichlet boundary condition is imposed at
$v_m$, and everywhere below we will be dealing with the reduced
boundary $\Gamma=\widetilde\Gamma\backslash\{ v_m\}.$ Since the
graph under consideration is a tree, for every $a,b\in\Omega,$
$a\not=b,$ there exists the unique path $\pi[a,b]$ connecting
these points.

For simplicity of the formulation of the balance conditions at the
internal vertexes, we introduce the special parametrization of
$\Omega$: we assume that at any internal vertex all the edges
connected at it have this vertex as start point or as end point.
We assume that clamped vertex $v_m$ is the start point of the edge
$e_m,$ which fix the parametrization.

Let $J:=\begin{pmatrix} 0&1\\-1&0\end{pmatrix}$, at each edge
$e_i$ we are given with a real matrix-valued potential
$V_i=\begin{pmatrix} p_i&q_i\\q_i&-p_i
\end{pmatrix}$, $p_i,q_i\in C^1(e_i)$. The space of real vector valued
square integrable functions on the graph $\Omega$ is denoted
by\linebreak $L_2(\Omega):= \bigoplus_{i=1}^{N}L_2(e_i,
\mathbb{R}^2).$ For the element $U \in L_2(\Omega)$ we write
\begin{equation*}
U:=\begin{pmatrix}u^1 \\u^2\end{pmatrix} =
        \left\{
            \left(
                \begin{matrix}
                    u^1_i \\
                    u^2_i
                \end{matrix}
            \right)
\right\}_{i=1}^N\ ,\ u^1_i, u^2_i \in L_2(e_i).
\end{equation*}
The continuity condition at the internal vertexes reads:
\begin{equation}
\label{Cont} u^1_i(v)=u^1_j(v), \quad \quad e_i\sim v,\,\, e_j\sim
v,\quad v\in V\backslash\widetilde\Gamma.
\end{equation}
%Here $e_{j} \sim v$ means edge $e_{j}$ is incident to vertex $v$,
The second condition (force balance) at the internal vertex $v$ is
introduced as
\begin{equation}
\label{Cont1} \sum_{i|e_i\sim v}u^2_i(v)=0, \quad v\in
V\backslash\widetilde\Gamma.
\end{equation}
We put $\Psi:=\begin{pmatrix}\psi^1\\ \psi^2\end{pmatrix} \in
L_2(\Omega) $,\, $\psi^1_i,\, \psi^2_i\in H^1(e_i)$ and introduce
the operator
\begin{equation*}
\mathcal{L}\Psi:=\left\{J\frac{d}{dx}\begin{pmatrix}\psi^1_i\\
\psi^2_i\end{pmatrix}+V_i\begin{pmatrix}\psi^1_i\\
\psi^2_i\end{pmatrix}\right\},\quad x\in e_i
\end{equation*}
with the domain
\begin{eqnarray*}
D(\mathcal{L})=\left\{\Psi\in L_2(\Omega)\Bigl|\, \psi^1_i,\,
\psi^2_i\in H^1(e_i),\,i=1,\ldots,N,\right.\\
\left.\Psi\,\, \mbox{satisfies}\ (\ref{Cont}),
(\ref{Cont1}),\,\,\psi^1(v)=0,\,v\in\Gamma \right\}
\end{eqnarray*}

By $\textbf{S}$ we denote the spectral problem on the graph:
\begin{align}
J\Psi_x+V\Psi=\lambda \Psi \ \ x \in e_i, \label{Dir_eqn} \\
\Psi\,\, \mbox{satisfies}\ \ (\ref{Cont}), (\ref{Cont1}) \
\mbox{at}
\ \        v\in V\backslash \widetilde\Gamma \label{Kirch}\\
\psi^1(v) = 0\ \ \text{for}\,\, v\in\Gamma , \label{Bound}
\end{align}
%It is known that the problem (\ref{Dir_eqn})--(\ref{Bound}) has a
%solution $\left\{\lambda_k,\Psi_k\right\}_{k=1}^\infty$, where
%$\Psi_k$ are chosen such that $\int_\Omega
%\Psi_k\Psi_l\,dx=\delta_{kl}$. We introduce the notation:
%$\gamma_k:=\left\{\gamma^k_1,\ldots\gamma^k_{m-1}\right\}=\left\{\psi_k^2(v_1)\ldots,\psi_k^2(v_{m-1})\right\}$.
%With (\ref{Dir_eqn})--(\ref{Bound}) we associate the
%\emph{spectral data}, the set of pairs:
%\begin{equation}
%\label{SD} \left\{\lambda_k,\gamma_k\right\}_{k=1}^\infty.
%\end{equation}
We introduce the Titchmarsh-Weyl (TW) matrix-function as an analog
to Dirichlet-to-Neumann map \cite{AK,AMN,ACLM} in the following
way: for $\lambda\notin \mathbb{R}$ and $\xi\in \mathbb{R}^{m-1}$
we consider the problem (\ref{Dir_eqn}), (\ref{Kirch}) with the
nonhomogeneous boundary condition:
\begin{equation}
\psi^1(v_i) = \xi_i,\quad i=1,\ldots,m-1.\label{Bound_nhm}
\end{equation}
The TW matrix function connects the values of the solution
$\Psi(\cdot,\lambda)$ to (\ref{Dir_eqn}), (\ref{Kirch}),
(\ref{Bound_nhm}) in the first and second channels at the
boundary:
\begin{eqnarray}
\label{TW_matr} \psi^2(\cdot,\lambda)|_{\Gamma}={\bf
M}(\lambda)\psi^1(\cdot,\lambda)|_{\Gamma},\\
\left(\psi^2(v_1,\lambda),\ldots,\psi^2(v_{m-1},\lambda)\right)^T
= {\bf M}(\lambda)\left(\xi_1,\ldots,\xi_{m-1}\right)^T.\notag
\end{eqnarray}
The inverse problem for the problem $\bf{S}$ is to recover the
tree $\Omega$, i.e. connectivity of edges and their lengths, and
parameters $p_i,$ $q_i$ on edges $e_i$ from $\bf{M}(\lambda)$.

Along with the spectral, we consider the dynamic inverse problem.
We introduce the \emph{outer space } $\mathcal{F}_\Gamma^T:=
L_2([0,T], \mathbb{R}^{m-1})$, the space of controls acting on the
reduced boundary of $\Omega$. The forward problem is described by
the Dirac system on the each edge of the tree:
\begin{equation}
iU_t(x,t)+JU_x(x,t)+V(x)U(x,t)=0 \ \ x \in e_i,\,\, t\geqslant
0,\label{Dir_eqn_time}
\end{equation}
conditions at internal vertexes:
\begin{equation}
U(v,t)\,\,\text{satisfies (\ref{Cont}), (\ref{Cont1}) for all
$t\geqslant 0$,\,\, $v\in
V\backslash\widetilde\Gamma,$}\label{Kirh_time}
\end{equation}
Dirichlet boundary conditions:
\begin{equation}
U^1|_{\Gamma} = F, \mbox{ on }\Gamma \times [0,T],\,\
\label{DirichletBoundCnd}
\end{equation}
where $F=\left(f^1(t),\ldots, f^{m-1}(t)\right)^T\in
\mathcal{F}^T_\Gamma$, and (\ref{DirichletBoundCnd}) means that
\begin{equation*}
\left(u^1(v_1,t),\ldots,
u^1(v_{m-1},t\right)^T=\left(f^1(t),\ldots, f^{m-1}(t)\right)^T.
\end{equation*}
By \textbf{D} we denote the dynamic problem on $\Omega$, described
by system (\ref{Dir_eqn_time}), compatibility conditions
\eqref{Kirh_time} at all internal vertices for all $t > 0$,
Dirichlet boundary condition \eqref{DirichletBoundCnd} and zero
initial condition $U(\cdot,0)=0$. The solution to this problem is
denoted by $U^F$. We introduce the \emph{response operator} for
the problem \textbf{D} by
\begin{equation}
\label{Resp_def} R^T \{F\}(t):= u^2(\cdot,t)|_\Gamma , \ \ t \in
[0,T].
\end{equation}
In other words, $R^T$ connects values of the solution $U^F$ to the
problem \textbf{D} in the first and the second channels at the
boundary: \begin{equation*} \left(R^T
\left(f^1(t),\ldots,f^{m-1}(t)\right)^T\right)(t)=\left(
u^2(v_1,t),\ldots, u^2(v_{m-1},t)\right)^T.
\end{equation*}
The operator $R^T$ has a form of convolution:
\begin{equation*}
\left(R^TF\right)(t)=\left(\mathbf{R}*F\right)(t)=\int_0^t
\mathbf{R}(t-s)F(s)\,ds,
\end{equation*}
where $\mathbf{R}(t)=\{R_{ij}\}_{i,j=1}^{m-1}$ is a \emph{response
matrix}. The entries $R_{ij}(t)$ are defined in the following way:
let $U_i$ be a solution to the boundary value problem
(\ref{Dir_eqn_time}), (\ref{Kirh_time}), $U_i(\cdot,0)=0$ with
special boundary condition (\ref{DirichletBoundCnd}) where
$F=(0,\ldots,\delta(t),\ldots,0)^T$ with only nonzero element at
$i-$th place. Then
\begin{equation}
\label{RespMD} R_{ij}(t)=u^2_i(v_j,t).
\end{equation}
The inverse problem for the problem $\bf{D}$ is to recover the
tree (connectivity of the edges and their lengths) and the matrix
potential on edges from the response operator $R^T(t)$, $t>0$
(\ref{Resp_def}).

The connection between spectral and the dynamic inverse data is
known \cite{AK,AMN,ACLM} and was used for solving inverse spectral
and dynamic problems. Let $F \in \mathcal{F}_\Gamma^T \cap
(C_0^{\infty}(0, +\infty))^{m-1}$ and
\begin{equation*}
\widehat{F}(k) := \int_{0}^{\infty} F(t) e^{ikt}dt
\end{equation*}
be its Fourier transform. The systems  (\ref{Dir_eqn})  and
(\ref{Dir_eqn_time}) are clearly connected: going formally in
(\ref{Dir_eqn_time}) over to the Fourier transform, we obtain
(\ref{Dir_eqn}) with $\lambda = k$. It is not difficult to check
that the response matrix function ${\bf R}(t)$ and TW matrix
function $\textbf{M}(\lambda)$ (Nevanlinna type matrix function)
are connected by the same transform:
\begin{equation}
\label{Sp_Dyn} \textbf{M}(k) = \int_{0}^{\infty} {\textbf R}(t)
e^{ikt}dt
\end{equation}
where this equality is understood in a weak sense. We use this
relationship between dynamic and spectral data solve the inverse
problem from either $\bf{M}(k)$ or $R(t)$, $t\geqslant 0.$

We will use the Boundary Control method \cite{BIP07}, first
applied to the problems on trees in \cite{B,BV} and its
modification, so called leaf-peeling method introduced in
\cite{AK} and developed in \cite{ACLM,AMN,2}. This method, as it
follows from its name is connected with the controllability
property of the dynamical system under the consideration. The
general principal \cite{BIP07} says that better controllability of
dynamical system leads to better identifiability. Introduce the
\emph{control operator}: $W^T: \mathcal{F}^T_\Gamma\mapsto
L_2(\Omega)$ acting by the rule:
\begin{equation*}
W^TF:=U^F(\cdot,T).
\end{equation*}
For the wave equation on a tree \cite{B,BV} the corresponding
control operator is boundedly invertible for certain values of
time. For the two-velocity system \cite{ACLM,ABCM} the
corresponding operator is not invertible, but at least there is
some "local" controllability. But for the Dirac system there are
no even "local" controllability, the latter causes the
consequences for the inverse problem. To overcome this difficulty,
we will use some ideas from \cite{BM_Dir}. In \cite{AMN} the
authors developed purely dynamic version of the leaf-peeling
method for the inverse problem for the wave equation with
potential on a finite tree, which allows one to solve the inverse
problem using $R^T$ for some finite $T.$ We are planning to return
to this (optimal in time) setting for a Dirac system on a tree
elsewhere.

The next section is devoted to the solution of IP. On analyzing
the reflection of a wave propagating from a boundary from an inner
vertex, we obtain the length of boundary edge. On the next step,
using the method from \cite{BM_Dir}, we find $p_i,$ $q_i$, for
boundary edges $e_i$ (i.e. $e_i\sim v_i$, $i=1,\ldots,m-1$). Then
we determine sheaf -- a star-shaped subgraph of $\Omega$,
consisting of boundary edges $e_1,\ldots, e_{m_0}$ and only one
non-boundary edge. And on the last step we consider the new tree:
$\Omega\backslash \cup_{i=1}^{m_0} e_i$ and recalculate the
Weyl-Titchmarsh matrix $\widetilde {\bf{M}}(\lambda)$ for this
reduced tree.

\section{Inverse problem.}

\subsection{Reflection from the inner vertex.}

Let $U^\delta$ be a solution to the special boundary value problem
for Dirac system on a tree: on the edges $e_i$, $i=1,\ldots,N$,
$U^\delta$ satisfies (\ref{Dir_eqn_time}), at internal vertexes
the continuity and force balance conditions (\ref{Kirh_time})
hold, $U^\delta(\cdot,T)=0$, and on the boundary we prescribe the
special condition
\begin{equation*}
u^1(v_1,t)=\delta(t),\,\,u^1(v_2,t)=0,\ldots, u^1(v_{m-1},t)=0.
\end{equation*}
We denote by $l_1$ the length of the boundary edge $e_1=[v_1,v']$,
which is identified with the interval $[0,l_1]$. We assume that
this edge is connected at the inner vertex $v'$ with other $n-1$
edges $e_2,\ldots,e_n$ which we identify with the intervals
$[l_1,l_1+l_i]$, where $l_i$ is a length of $e_i$, $i=2,\ldots,n$.
When $t<l_1$ (i.e. the wave generated at $v_1$ does not each the
inner vertex $v'$), the solution to the above problem is zero on
all edges except $e_1$. And on this edge it is given by
\cite{BM_Dir}:
\begin{equation*}
U_1^\delta(x,t)=\begin{pmatrix} \delta(t-x) \\
i\delta(t-x)\end{pmatrix}+\Gamma(x,t),\quad x\in e_1,\,\, t<l_1,
\end{equation*}
where $\Gamma(x,t)$ is a smooth function in the region
$\{0<x<t\}$. At $t=l_1$ the wave reaches the inner vertex $v'$, on
the time interval $l_1<t<l_1+L$ where
$L=\min_{i=2,\ldots,n}\{l_i\}$, the solution on the edges
$e_2,\ldots, e_n$ has a form
\begin{equation*}
U_i^\delta(x,t)=\begin{pmatrix} \alpha\delta(t+x-l_1-l_i) \\
-i\alpha\delta(t+x-l_1-l_i)\end{pmatrix}+\Gamma_i(x,t),\quad x\in
e_i,\,\, l_1<t<l_1+L.
\end{equation*}
and on the first edge, $x\in e_1$:
\begin{equation*}
U_1^\delta(x,t)=\begin{pmatrix} \delta(t-x) \\
i\delta(t-x)\end{pmatrix}+\begin{pmatrix} \gamma\delta(t+x-2l_1) \\
-i\gamma\delta(t+x-2l_1)\end{pmatrix}+\Gamma'(x,t),\ \,\,
l_1<t<l_1+L.
\end{equation*}
In the above representations $\Gamma_i$, $\Gamma'$, $i=2,\ldots,n$
are smooth functions and the constants $\alpha, \gamma$ are
subjected to determination.

We use the first continuity condition (\ref{Cont}) at the vertex
$v'$ to get the relation
\begin{equation*}
1+\gamma=\alpha,
\end{equation*}
and use the force balance condition (\ref{Cont1}) at $v'$, which
yields
\begin{equation*}
1-\gamma-(n-1)\alpha=0.
\end{equation*}
Two above equalities lead to the following formulas
\begin{equation*}
\gamma=\frac{2-n}{n},\quad \alpha=\frac{2}{n}.
\end{equation*}
Bearing in mind the definition of a response matrix
(\ref{RespMD}), we see that its component $R_{11}(t)$ has a form:
\begin{equation}
\label{UDelta_rsp} R_{11}(t)= u^\delta_2(v_1,t)=
i\delta(t)+i\frac{n-2}{n}\delta(t-2l_1)+\Gamma'_2(0,t),\quad t\in
[0,2l_1+L)
\end{equation}
with some smooth $\Gamma'_2(0,t)$. Thus by knowing the response
matrix entry $R_{11}(t)$, one can determine the length $l_1$ of
the edge $e_1$ (it is contained in the argument of the second
singular term) and the number of edges, $e_1$ connected with. The
representation (\ref{UDelta_rsp}) implies that from the diagonal
elements $R_{ii}$ $i=1,\ldots,m-1$ of the response matrix one can
extract the lengths $l_i$ of the boundary edges $e_i$,
$i=1,\ldots,m-1$.

\subsection{Inverse problem on a half-line.}

Here we show following \cite{BM_Dir} that diagonal elements of the
response matrix determine not only the lengths $l_i$ of boundary
edges $e_i$, but also matrix potentials $V_i$ on $e_i$,
$i=1,\ldots, m-1$.

We consider the inverse problem for the Dirac system on a
half-line, which is set up in the following way:
\begin{equation}
\label{Dir_int}
\left\{
\begin{array}l
iu_t+Ju_x+Vu=0, \quad  x>0, \quad 0<t<2T \\
u\big|_{t=0}=0, \quad  x \geqslant 0 \\
u_1\big|_{x=0}=f, \quad 0\leqslant t \leqslant 2T\,,
\end{array}
\right.
\end{equation}
where $V=\begin{pmatrix} p&q\\q&-p
\end{pmatrix}$ is a matrix potential, $p=p(x)$ and $q=q(x)$ being
real-valued $C^1$-smooth functions. We associate a {\it response
operator} to the above system, acting in $L_2([0,2T]; {\mathbb
C})$ by the rule
\begin{equation*}
R f:=u^f_2(0,\cdot)\big|_{0 {\leqslant} t {\leqslant} 2T}.
\end{equation*}
This operator has a form of convolution: $Rf=if+r\ast f$, where
$r\big|_{0 {\leqslant} t {\leqslant} 2T}$ is a {\it response
function}. The response function $r(t)$ for $t\in (0,2T)$ is
determined by the values of the potential $V(x)$ for $x\in (0,T)$
only, therefore, the relevant dynamic setup of the inverse problem
is: for a fixed $T>0$, given $r\big|_{0 {\leqslant} t {\leqslant}
2T}$, to recover $V\big|_{0{\leqslant} x {\leqslant} T}$. We
assume that all functions of time $t\geqslant 0$ are extended to
$t<0$ by zero. Also, for a $z \in \mathbb C$, $\bar
z:=\operatorname{Re} z -i \operatorname{Im} z$ is its conjugate.

A $\mathbb C$-valued function $r\big|_{0 {\leqslant} t {\leqslant}
2T}$ determines an operator $C^T$ acting in\linebreak $L_2([0,2T];
{\mathbb C}^2)$ by the rule
\begin{equation}
\label{C_T} \left(C^Ta\right)(t)\,=\,2 a(t)+\int_0^T c^T(t,s)
a(s)\,ds\,, \qquad 0\leqslant t \leqslant T,
\end{equation}
where $a=\begin{pmatrix}a_1(t)\\a_2(t)\end{pmatrix}$, and the the
element of  matrix kernel $c^T$ are
 \begin{align*}
\notag & c_{11}(t,s)(t)=-i\,[r(t-s)-\bar r(s-t)]\,, \quad
c_{12}(t,s)=-i\,\bar r(2T-t-s)\,,\\
& c_{21}(t,s)= i\,r(2T-t-s)\,,\quad c_{22}(t,s)(t)=i\,[\bar
r(t-s)-r(s-t)]\,.
 \end{align*}
In \cite{BM_Dir} it is proved
\begin{theorem}
The function $r\in C^1([0,2T]; {\mathbb C})$ is the response
function of a system (\ref{Dir_int}) with a {\it $C^1$-smooth real
zero trace} potential $V$ if and only if operator $C^T$ is a
positive definite isomorphism.
\end{theorem}
Below we describe a procedure of recovering a potential $V$ from a
response function $r$:

\noindent{\bf$1.$}\,\,\,Given a response function
$r(t),\,\,0\leqslant t \leqslant 2T$ of the system
(\ref{Dir_int}), determine the operator $C^T$ and the
matrix-kernel $c^T$ by (\ref{C_T}).
\smallskip

\noindent{\bf$2.$}\,\,\, For $0<\xi<T$ solve the family of the
{\it linear integral equations}
\begin{equation}
\label{GL_eqn}
\frac{1}{2} c^T(t,s)+2k^\xi(t,s)+\int_{T-\xi}^T
k^\xi(t,\eta) c^T(\eta,s)\,d\eta=0\,, \quad T-\xi \leqslant s,t
\leqslant T\,,
\end{equation}
which determines the matrix-valued function
$k^\xi,\,\,\,0<\xi\leqslant T$ via $c$. The solvability is
guaranteed by the positive-definiteness of $C^T$. By standard
integral equations theory arguments, the solution $k^\xi$ is of
the same smoothness as $c^T$, i.e., is $C^1$-smooth outside the
diagonal $t=s$.
\smallskip

\noindent{\bf$3.$}\,\,\,define the matrix $ w(x,x)$ by
\begin{equation}
\label{diag check w}  w(x,x)=\, - 2\begin{pmatrix} 1 & 1\\i &
-i\end{pmatrix} k^x(T-x, T-x)\,, \qquad x \in (0,T),
\end{equation}
take its first column $\begin{pmatrix} w_1 \\ w_2\end{pmatrix}$
and recover the entries of matrix potential by:
\begin{eqnarray*}
p(x)&=\operatorname{Im} w_1(x,x)+ \operatorname{Re} w_2(x,x),\quad
x \in (0,T),\\
q(x)&=-\operatorname{Re} w_1(x,x)+ \operatorname{Im}
w_2(x,x),\quad x \in (0,T).
\end{eqnarray*}
We use the method described above to recover the potential $V_i $
on each boundary edge $e_i$, $i=1,\ldots,m-1$. Where for a fixed
boundary vertex $v_i$, we consider the response $R^{2T_i}_{ii}$
with $T_i=l_i$ and $l_i$ is a length of $e_i$, were recovered from
$R_{ii}(t)$ as explained in the previous subsection.

\subsection{Recovery of the boundary sheaves.}

At this point we assume that we already know the lengths $l_i$ of
boundary edges $e_i$, $i=1,\ldots,m-1$.  We will need the reduced
response function ${\bf R}(t)=\{R_{ij}(t)\}_{i,j=1}^{m-1}$
(\ref{Resp_def}), (\ref{RespMD}). If our inverse data is a TW
function $\bf{M}(\lambda)$, we can pass to $R(t)$ taking the
inverse Fourier transform of (\ref{Sp_Dyn}).

First of all we identify the edges connected at the same vertex.
Two boundary edges, $e_i$ and $e_j,$ $1\leqslant i,j\leqslant m-1$
have a common vertex if and only if
\begin{equation*}
R_{ij}(t)=\left\{\begin{array}l =0,\quad \text{for}\,\, t<l_i+l_j,\\
\not=0,\quad \text{for}\,\, t>l_i+l_j.
\end{array}\right.
\end{equation*}
This relation allows us to divide the boundary edges into groups,
such that edges from one group have a common vertex. We call these
groups pre-sheaves.  More exactly, we introduce the following

\begin{definition}  We consider a subgraph of $\Omega$ which
is a star graph consisting of {\sl all} edges incident to an
internal vertex $v.$ This star graph is called a \emph{pre-sheaf}
if it contains at least one  boundary edge of  $\Omega.$ A
pre-sheaf is called a \emph{sheaf} if all but one its edges are
the boundary edges of  $\Omega.$
\end{definition}

The sheaves are especially important to our identification
algorithm. To extract them we denote the found pre-sheaves by
$P_1, \ldots, P_L$, and define  the distance $d(P_k,P_m)$ between
two pre-sheaves in the following way: we take boundary edges
$e_i\in P_k$ and $e_j\in P_m$ and then  put
\begin{equation*}
d(P_k,P_m)=\max\{t>0\, : \, R_{ij}(t-l_i-l_j)=0\}.
\end{equation*}
Clearly this definition does not depend on the particular choice
of $e_i\in P_k$ and $e_j\in P_m$ and gives the distance between
the internal vertices of the
 pre-sheaves
$P_k $ and $ P_m.$ Then we consider
\begin{equation*}
\max_{k,m\in 1,\ldots,N,\,k\not=m}d(P_k,P_m).
\end{equation*}
It is not difficult to see that two pre-sheaves on which this
maximum is attained (we denote them by $P$ and $P'$) are sheaves.
Indeed, since $\Omega$ is a tree, there is only one path between
$P$ and $P'$. If we assume the existence of  an ``extra" internal
edge in $P$ or $P'$, this leads to contradiction, since there
would necessarily exist sheaves with a distance between them
greater than $d(P,P')$.

\subsection{Leaf peeling method. }

Let the sheaf $P$, found on the previous step, consist of the
boundary vertices $v_1,\ldots,v_{m_0}$ from $\Gamma$, the
corresponding boundary edges $e_1,\ldots,e_{m_0}$ and an internal
edge $e_{m'_0}.$ We assume that we already recovered lengths $l_i$
and potentials $V_i$, $i=1,\ldots,m'_0$, i.e. on boundary edges of
$P$. We identify each edge $e_{m'_0},e_i, \, i=1,\ldots,m_0$, with
the interval $[0,l_i]$ and the vertex $v_{m'_0},$  the internal
vertex of the sheaf, ---  with the set of common endpoints $x=0.$
At this point it is convenient to renumerate the edge $e_{m'_0}$
as $e_{0}$ and the vertex $v_{m'_0}$ as  $v_{0}$.

By $\widetilde {\mathbf{M}}(\lambda)$ we denote the reduced TW
matrix function associated with the new graph $\widetilde
\Omega=\Omega\backslash \bigcup_{i=1}^{m_0}\{e_i\}$  with boundary
points $v_0\cup\Gamma\backslash \bigcup_{i=0}^{m_0}v_i$.

First we recalculate entries $\widetilde M_{0i}(\lambda)$,
$i=0,m_0+1,\ldots,m-1$. Let us fix $v_1$, the boundary point of
the sheaf $P$. Let $\Psi$ be a solution to the problem
(\ref{Dir_eqn})-(\ref{Kirch}) with the boundary conditions given
by
\begin{equation*}
\psi^1(v_1)=1,\quad \psi^1(v_j)=0,\quad j=2,\ldots,m.
\end{equation*}
We point out that on the boundary edge $e_1$ the function $\Psi$
solves the Cauchy problem
\begin{equation}
\label{v1}
\left\{\begin{array}l
J\Psi_x+V\Psi=\lambda \Psi \quad x\in e_1,\\
\psi^1(v_1)=1,\,\, \psi^2(v_1)=M_{11}(\lambda).
\end{array}
\right.
\end{equation}
On other boundary edges of $P$, the function $\Psi$ solves the
problems
\begin{equation}
\label{v3} \left\{
\begin{array}l
J\Psi_x+V\Psi=\lambda \Psi,\quad x\in e_i, \\
\psi^1(v_i)=0,\,\, \psi^2(v_i)=M_{1i}(\lambda), \quad
i=2,\ldots,m_0
\end{array}
\right.
\end{equation}
Since we know potentials $V_i$ on the edges $e_1,\ldots,e_{m_0}$,
we can solve the Cauchy problems (\ref{v1}) and (\ref{v3}), and
use the conditions (\ref{Cont}), (\ref{Cont1}) at the internal
vertex $v_{0}$ to recover $\psi^1_0(v_0,\lambda)$,
$\psi^2_0(v_0,\lambda)$ -- the value of the solution $\Psi$ at
$v_0$, i.e. at the ``new" boundary point of the new tree
$\widetilde\Omega$. Then we obtain:
\begin{eqnarray*}
    &\widetilde M_{00}(\lambda)=\frac{\psi^2_0(v_0,\lambda)}{\psi^1_0(v_0,\lambda)},\\
    &\widetilde M_{0i}(\lambda)=\frac{M_{1i}(\lambda)}{\psi^1_0(v_0,\lambda)},\quad
    i=m_0+1,\ldots m-1.
\end{eqnarray*}
To find $\widetilde M_{i0}(\lambda)$, $i=m_0+1,\ldots,m-1$ we fix
boundary point $v_i$, $i\notin\{1,\ldots,m_0,m\}$ and consider the
solution $\Psi$ to (\ref{Dir_eqn})--(\ref{Kirch}) with the
boundary conditions given by
\begin{equation*}
\psi^1(v_i)=1,\quad \psi^1(v_j)=0,\quad j\not=i.
\end{equation*}
The function $\Psi$ solves the Cauchy problems on the edges
$e_1,\ldots,e_{m_0}$:
\begin{equation}
\label{v5} \left\{
\begin{array}l
J\Psi_x+V\Psi=\lambda \Psi,\quad x\in e_j\\
\psi^1(v_j)=0,\,\, \psi^2(v_j)=M_{ij}(\lambda),\quad
j=1,\ldots,m_0.
\end{array}
\right.
\end{equation}
Since we know the potential on the boundary edges of $P$, we can
solve Cauchy problems (\ref{v5}) and use conditions at the
internal vertex $v_0$ to recover $\psi^1_0(v_0,\lambda)$,
$\psi^2_0(v_0,\lambda)$ -- the value of solution at the ``new''
boundary point $v_0$ of reduced tree $\widetilde \Omega$.

On the the other hand, on new tree $\widetilde\Omega$ the function
$U$ solves the problem (\ref{Dir_eqn})-(\ref{Kirch}) with the
boundary conditions
\begin{eqnarray*}
\psi^1(v_i)=1,\,\, \psi^1(v_{0})=\psi_0^1(v_0,\lambda), \,\,
\psi^1(v_j)=0,\, j=m_0+1,\ldots, m,\\
v_j\not=v_i,\,v_j\not=v_{0}.
\end{eqnarray*}
Thus for the entries of $\widetilde{{\bf M}}(\lambda)$ the
following relations hold:
\begin{eqnarray*}
\widetilde
M_{i0}(\lambda)=\psi_0^2(v_0,\lambda)-\psi_0^1(v_0,\lambda)
\widetilde M_{00}(\lambda),\\
\widetilde M_{ij}(\lambda)=M_{ij}(\lambda)-\psi_0^1(v_0,\lambda)
\widetilde M_{0j}(\lambda).
\end{eqnarray*}
To recover all elements of the reduced matrix $\widetilde{{\bf
M}}(\lambda)$ we need to repeat this procedure for all
$i,j=m_0+1,\ldots,m-1$.

Thus using the  described procedure we can recalculate the
truncated TW matrix $\widetilde{{\bf M}}(\lambda)$ for the new
`peeled' tree $\widetilde\Omega$. Repeating the procedure
sufficient number of times, we step by step recover the tree and
the matrix potential.

\noindent{\bf Acknowledgments}

The research of Victor Mikhaylov was supported in part by RFBR
17-01-00529. Alexandr Mikhaylov was supported by RFBR 17-01-00099;
A. S. Mikhaylov and V. S. Mikhaylov were partly supported by VW
Foundation program ``Modeling, Analysis, and Approximation Theory
toward application in tomography and inverse problems.'' Gulden
Murzabekova and Victor Mikhaylov were also partly supported by the
Ministry of Education and Science of Republic of Kazakhstan, grant
no. 4290/GF4.

\end{document}